\documentclass[11pt,a4paper]{article}
\usepackage{amsmath, amsfonts, amssymb,amsthm}
\usepackage{graphicx}
\usepackage{hyperref}
\usepackage[utf8]{inputenc}  

\newtheorem{thm}{Theorem}

\newtheorem{rem}[thm]{Remark}
\newtheorem{conj}[thm]{Conjecture}

\title{Pascal Determinantal Arrays and A Generalization of Rahimpour's Determinantal Identity}
\author{H. Teimoori \and H. Khodakarami}
\date{}

\begin{document}

\maketitle

\begin{abstract}
We introduce a new infinite family of arrays, the \emph{Pascal determinantal arrays} of order $k$, denoted $PD_k$, which generalize the classical Pascal array via determinantal constructions. We present a recursive algorithm for generating $PD_k$, establish its correctness using Dodgson's condensation and a weighted sliding-cross rule, and provide a geometric interpretation of the entries $P^{(k)}_{i,j}$ as weighted double sticks in the Pascal plane. Our main result proves a conjecture that generalizes Rahimpour's identity: for all $i,j,k \ge 0$,
\[
P^{(k)}_{i,j} = P^{(j)}_{i,k},
\]
where $P^{(k)}_{i,j}$ is the determinant of the $k \times k$ subarray of the Pascal array starting at $(i,j)$. The proof combines algebraic recurrence techniques with a visual, geometry-based argument that reveals the intrinsic symmetry of determinantal Pascal structures.
\end{abstract}


\section{Introduction}

Pascal's triangle is a cornerstone of combinatorics, rich with algebraic, geometric, and number-theoretic patterns. Among its many fascinating properties, determinantal identities within the triangle have attracted considerable attention. Rahimpour's conjecture---that the determinant of any square subarray with its left edge in the first column equals its top-right entry---was proved by Teimoori and Bayat \cite{m1} using matrix-theoretic methods. This result naturally invites a broader question: what happens when the left edge lies in an arbitrary column $k$?

In this paper, we answer this question by introducing \emph{Pascal determinantal arrays} $PD_k$, defined recursively as arrays whose entries are determinants of $k \times k$ subarrays of the classical Pascal array. We prove that these arrays satisfy an elegant symmetric identity that extends Rahimpour's observation to all $k \ge 1$:

\begin{equation}\label{eq:main}
P^{(k)}_{i,j} = P^{(j)}_{i,k} \qquad (i,j,k \ge 0).
\end{equation}

Our approach is multifaceted. In Section~\ref{sec:algorithm}, we give a simple recursive algorithm that generates $PD_k$ from $PD_{k-1}$. Section~\ref{sec:dodgson} recalls Dodgson's condensation formula and derives a recurrence for $PD_k$ entries. Section~\ref{sec:sliding} introduces a weighted version of the star-of-David rule, which we extend to all $PD_k$ via induction. Section~\ref{sec:correctness} proves the correctness of the recursive algorithm using the sliding-cross property. Finally, in Section~\ref{sec:geometry}, we provide a geometric interpretation of $P^{(k)}_{i,j}$ as the weight of a \emph{double stick} in the Pascal plane---a viewpoint that leads to a purely visual proof of \eqref{eq:main}.

The paper thus unites algebraic recurrence, combinatorial condensation, and geometric visualization to fully resolve the generalized Rahimpour conjecture, while uncovering new structural layers within Pascal's classical triangle.



\section{The Class of Pascal Determinantal Array}\label{sec:class}

We begin with the classical Pascal infinite \emph{square} array 
\[
P = \big( P_{i,j} = \binom{i+j}{i} \big)_{i,j \geq 0}.
\]
From $P$ we define a new family of infinite arrays, the \emph{Pascal determinantal arrays of order $k$}, denoted
\[
PD_k = \big( P^{(k)}_{i,j} \big)_{i,j \geq 0},
\]
where $P^{(k)}_{i,j}$ is the determinant of the $k \times k$ sub‑array of $P$ with top‑left corner at $(i,j)$:
\[
P^{(k)}_{i,j}:=
\left|
\begin{array}{ccc}
P_{i,j}       & \cdots & P_{i,j+k-1}   \\
\vdots        & \ddots & \vdots        \\
P_{i+k-1,j}   & \cdots & P_{i+k-1,j+k-1}
\end{array}
\right|.
\]

\begin{figure}[ht]
	\centering
	\includegraphics[scale=.75]{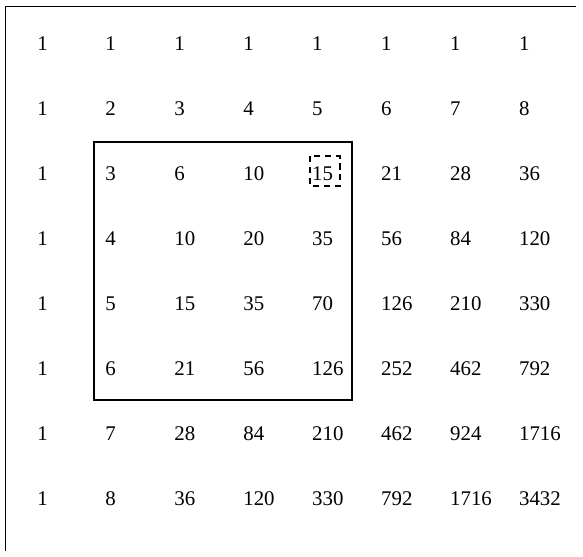}
	\caption{The Pascal array $P = PD_1$.}
	\label{fig:F1}
\end{figure}

\begin{figure}[ht]
	\centering
	\includegraphics[scale=.75]{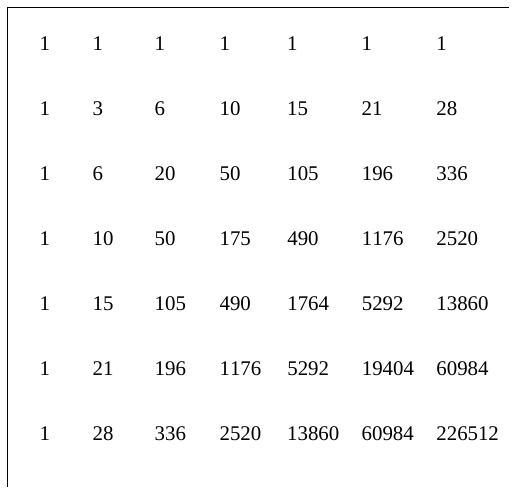}
	\caption{Pascal determinantal array $PD_2$ (Narayana triangle in square form).}
	\label{fig:F2}
\end{figure}

By definition, $PD_1$ coincides with the Pascal array itself (Figure~\ref{fig:F1}).  
The array $PD_2$, shown in Figure~\ref{fig:F2}, is known to be the square‑form of the Narayana triangle (sequence A001263 in the OEIS \cite{m5}).

\subsection{From Rahimpour's Identity to a Natural Generalization}

Rahimpour's determinantal identity concerns square sub‑arrays of $P$ whose left edge lies in the first column (Figure~\ref{fig:F1}).  
It states that the determinant of such an $n\times n$ sub‑array equals its top‑right entry.

\begin{thm}[Rahimpour's identity \cite{m1}]\label{RahimpourDetIdent}
	For the Pascal array $P = (P_{i,j})$,
	\[
	P^{(1)}_{i,j} = P^{(j)}_{i,1} \qquad (i,j \ge 0).
	\]
\end{thm}

A natural question arises: what happens when the left edge is placed in an arbitrary column $k$?  
While studying Rahimpour's result, the second author conjectured the following symmetric generalization.

\begin{figure}[ht]
	\centering
	\includegraphics[scale=.75]{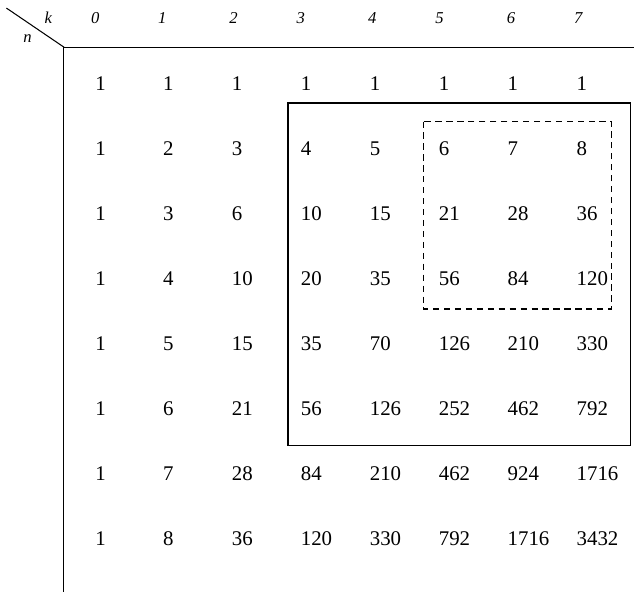}
	\caption{Illustration of Conjecture~\ref{HasanConj}: the determinant of the large square equals that of the small shaded square.}
	\label{fig:F3}
\end{figure}

\begin{conj}[Generalized Rahimpour identity]\label{HasanConj}
	For all $i,j,k \ge 0$,
	\[
	P^{(k)}_{i,j} = P^{(j)}_{i,k}.
	\tag{2}
	\]
	Equivalently, the determinant of an $n\times n$ sub‑array with left edge in column $k$ equals the determinant of the $k\times k$ sub‑array at its top‑right corner (Figure~\ref{fig:F3}).
\end{conj}

Theorem~\ref{RahimpourDetIdent} was proved in \cite{m1} using linear‑algebraic properties of Pascal functional matrices together with \emph{Jacobi adjoint formula}. 
In contrast, our proof of Conjecture~\ref{HasanConj} (given in Section~\ref{sec:geometry}) employs a geometric approach that relies on the weighted sliding‑cross rule (Section~\ref{sec:sliding}) and Dodgson's condensation (Section~\ref{sec:dodgson}), leading to a visual interpretation of $P^{(k)}_{i,j}$ as a weighted double stick in the Pascal plane.



\section{A Recursive Algorithm for Pascal Determinantal Arrays}\label{sec:algorithm}

In this section, we present a simple recursive algorithm that generates Pascal determinantal arrays of any order \(k\) directly from the standard Pascal array. The algorithm proceeds as follows:

\begin{enumerate}
	\item Set \(PD_1 :=\) The Pascal Array.
	\item Remove the zeroth row and zeroth column of \(PD_{k-1}\) and denote the resulting infinite array by \(RD_k = \big(R_{i,j}^{(k)}\big)_{i,j \geq 0}\).
	\item Define \(QD_{i,j}^{(k)} := \dfrac{R_{i,j}^{(k)}}{R_{0,i+j}^{(k)}}\) for all \(i,j \geq 0\).
	\item Set \(P_{i,j}^{(k)} := QD_{i,j}^{(k)} \cdot P_{i,j}^{(1)}\) for all \(i,j \geq 0\).
\end{enumerate}

The resulting array \(PD_k = \big(P_{i,j}^{(k)}\big)_{i,j \geq 0}\) is the Pascal determinantal array of order \(k\).

\medskip

Equivalently, the entries of \(PD_k\) satisfy the following recursive definition:

\[
PD_k = \big( P^{(k)}_{i,j} \big)_{i,j \geq 0},
\]
where
\begin{enumerate}
	\item \(P^{(1)}_{i,j} = \dbinom{i+j}{i}\),
	\item For \(k \geq 1\),
	\[
	P^{(k+1)}_{i,j} := \dfrac{P^{(k)}_{i+1,j+1} \; P^{(1)}_{i,j}}{P^{(k)}_{1,i+j+1}}.
	\]
\end{enumerate}

The two formulations are equivalent; the latter provides a compact recurrence that will be used extensively in subsequent sections.

\subsection*{Illustrative Example: Constructing \(PD_3\)}

The following steps demonstrate the construction of \(PD_3\) from \(PD_2\):

\begin{enumerate}
	\item[\textbf{Step 1.}] Start with \(PD_2\) (the Narayana array):
	
	\centering
	\includegraphics[scale=.8]{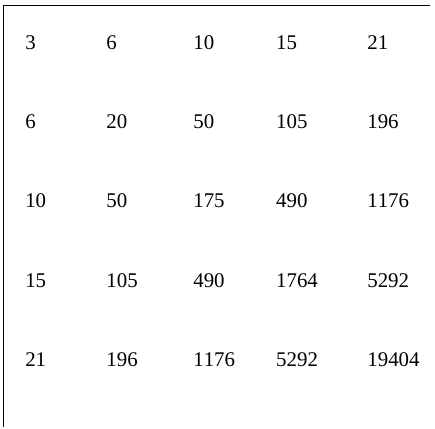}
	
	\item[\textbf{Step 2.}] Remove the zeroth row and column to form \(RD_3\):
	
	\centering
	\includegraphics[scale=.8]{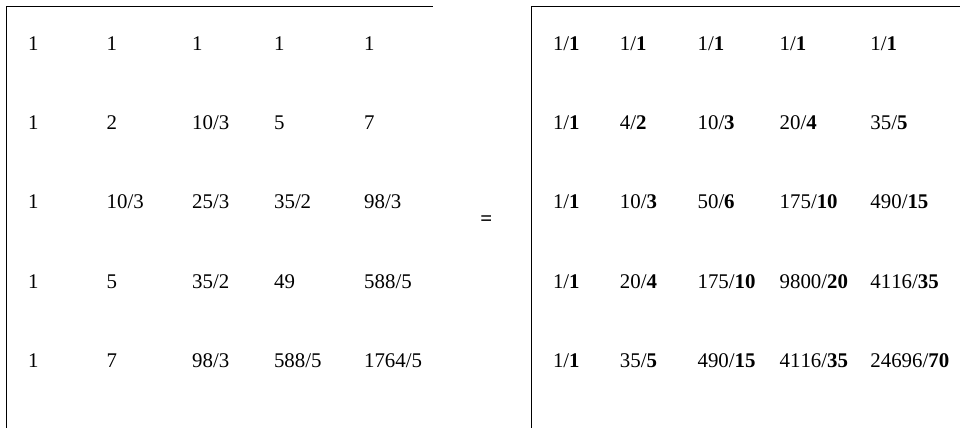}
	
	\item[\textbf{Step 3.}] Normalize and multiply by \(PD_1\) to obtain \(PD_3\):
	
	\centering
	\includegraphics[scale=.8]{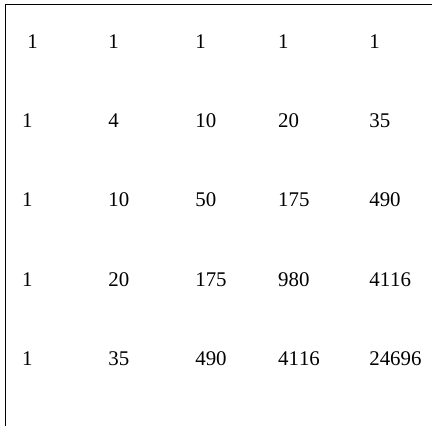}
	
\end{enumerate}

The correctness of this algorithm is not immediately obvious; it will be established in Section~\ref{sec:correctness} using a combination of algebraic and geometric properties of Pascal determinantal arrays, including Dodgson's condensation (Section~\ref{sec:dodgson}) and the weighted sliding‑cross rule (Section~\ref{sec:sliding}).



\section{Dodgson's Condensation of Determinants}\label{sec:dodgson}

Dodgson's condensation, introduced by Charles Dodgson in 1866 \cite{m4}, is a recursive method for computing determinants. For a square matrix \(A\) of order \(k \ge 3\), the condensation formula states:

\[
A_{k-2}(2,2) \, A_{k}(1,1) = A_{k-1}(1,1) \, A_{k-1}(2,2) - A_{k-1}(1,2) \, A_{k-1}(2,1),
\]

where \(A_r(i,j)\) denotes the \(r \times r\) minor formed from \(r\) contiguous rows and columns of \(A\), starting at row \(i\) and column \(j\). In particular:
\begin{itemize}
	\item \(A_k(1,1) = \det(A)\) (the full determinant),
	\item \(A_{k-2}(2,2)\) is the \emph{central minor},
	\item \(A_{k-1}(1,1)\), \(A_{k-1}(2,2)\), \(A_{k-1}(1,2)\), and \(A_{k-1}(2,1)\) are the northwest, southeast, northeast, and southwest minors, respectively.
\end{itemize}

\subsection*{Illustrative Example}

Consider a \(4 \times 4\) matrix
\[
A = 
\begin{bmatrix}
a_1 & a_2 & a_3 & a_4 \\
b_1 & b_2 & b_3 & b_4 \\
c_1 & c_2 & c_3 & c_4 \\
d_1 & d_2 & d_3 & d_4 
\end{bmatrix}.
\]
The relevant minors are:
\begin{align*}
A_2(2,2) &= 
\begin{bmatrix}
b_2 & b_3 \\
c_2 & c_3 
\end{bmatrix}, &
A_3(1,1) &=
\begin{bmatrix}
a_1 & a_2 & a_3 \\
b_1 & b_2 & b_3 \\
c_1 & c_2 & c_3
\end{bmatrix},\\[4pt]
A_3(2,2) &=
\begin{bmatrix}
b_2 & b_3 & b_4 \\
c_2 & c_3 & c_4 \\
d_2 & d_3 & d_4 
\end{bmatrix}, &
A_3(1,2) &=
\begin{bmatrix}
a_2 & a_3 & a_4 \\
b_2 & b_3 & b_4 \\
c_2 & c_3 & c_4 
\end{bmatrix},\\[4pt]
A_3(2,1) &=
\begin{bmatrix}
b_1 & b_2 & b_3 \\
c_1 & c_2 & c_3 \\
d_1 & d_2 & d_3
\end{bmatrix}.
\end{align*}

Thus every \(k \times k\) minor can be expressed in terms of \((k-1) \times (k-1)\) and \((k-2) \times (k-2)\) minors.

\subsection*{Recurrence for Pascal Determinantal Arrays}

Applying Dodgson's condensation to the Pascal array yields a recurrence relation for the entries of \(PD_k\). For all \(i,j \ge 0\) and \(k \ge 1\),

\begin{equation}\label{eq:dodgson-recurrence}
P_{i,j}^{(k)} = \frac{P_{i+1,j+1}^{(k-1)} \, P_{i,j}^{(k-1)} - P_{i+1,j}^{(k-1)} \, P_{i,j+1}^{(k-1)}}{P_{i+1,j+1}^{(k-2)}},
\end{equation}

with the initial arrays defined as:
\[
PD_0 := J \quad (\text{the all‑ones matrix}), \qquad PD_1 := P \quad (\text{the Pascal array}).
\]

This recurrence will play a central role in the proof of the algorithm's correctness (Section~\ref{sec:correctness}) and in the geometric interpretation of Section~\ref{sec:geometry}.


\section{The Weighted Version of Star of David Rule}\label{sec:sliding}

\begin{figure}[h!]
	\begin{center}
		\includegraphics[scale=0.75]{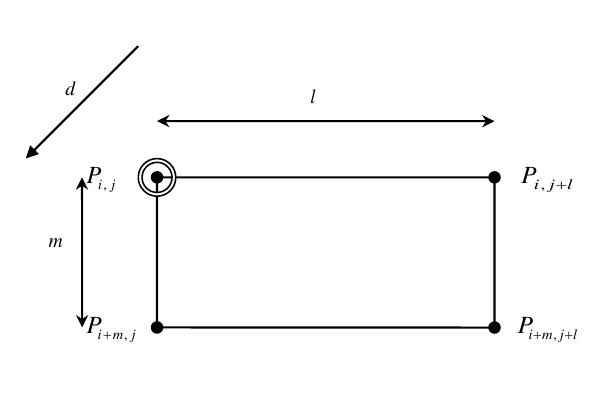}
		\caption{A rectangle with the Pascal array entries}
		\label{fig:F4}
	\end{center}
\end{figure}

Consider the Pascal array. We can draw an arbitrary rectangle whose vertices are entries of the Pascal array (see Figure~\ref{fig:F4}). We identify a vertex (the circled vertex) as the anchor of this rectangle. 
Now we define a weight for this rectangle denoted by
$W$, as follows

\begin{center}
	$W:=\frac{P_{i+m,j+l}.P_{i,j}}{P_{i+m,j}.
		P_{i,j+l}},~~~ m,l\geq 1,~~~ i,j\geq 0.$
\end{center}

Hilton and Pedersen proved \cite{m4} that when we move the anchor of the rectangle through the diagonal of the Pascal array (indicated by the arrow $d$ in Figure~\ref{fig:F4}), the weight remains constant. We will call
this property as the weighted version of star of David rule. 

\begin{figure}[h!]
	\begin{center}
		\includegraphics[scale=0.5]{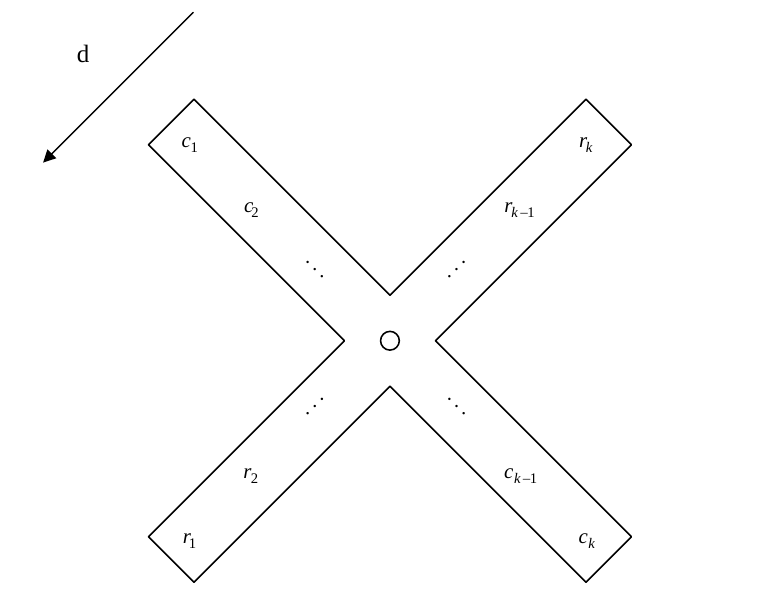}
		\caption{A cross of size $k$}
		\label{fig:F5}
	\end{center}
\end{figure}

We can define a symmetric cross of size $k$ in the Pascal array (See Figure~\ref{fig:F5}) as an $2k$-tuples ($c_1,\ldots,c_k,r_1,\ldots,r_k$), where 
$c_i,r_i,~(1\leq i \leq k)$ are the entries of the Pascal array. We then define the weight of such a cross, $W_c$, as

\begin{center}
	$W_c:=\frac{c_1.c_2. \ldots .c_k}{r_1.r_2. \ldots .r_k}.$
\end{center}

\begin{rem}
	
	An immediate consequence of the weighted version of star of David rule in the Pascal array is that if we move the center of a cross through the diagonal $d$, then the weight $W_c$ remains constant. We will call this property the \emph{sliding-cross rule} in the Pascal array.
	We also note that if the sliding-cross rule holds for $m=l=1$, then by an straight forward induction on $m$ and $l$, it can be easily seen that it holds in general for every positive integers $m$ and $l$.
	
\end{rem}

Our next step is to show that the above sliding rule is also true
for any Pascal determinantal array $PD_{k}$ $(k \geq 1)$. To do this,
we apply mathematical induction on $k$. The basis step, $k=1$, is
already true. For the inductive step, we assume that
every Pascal determinantal array of order less than $k+1$ has the
sliding property and we show that it is also true for $PD_{k+1}$.

\begin{figure}[h!]
	\begin{center}
		\includegraphics[scale=0.75]{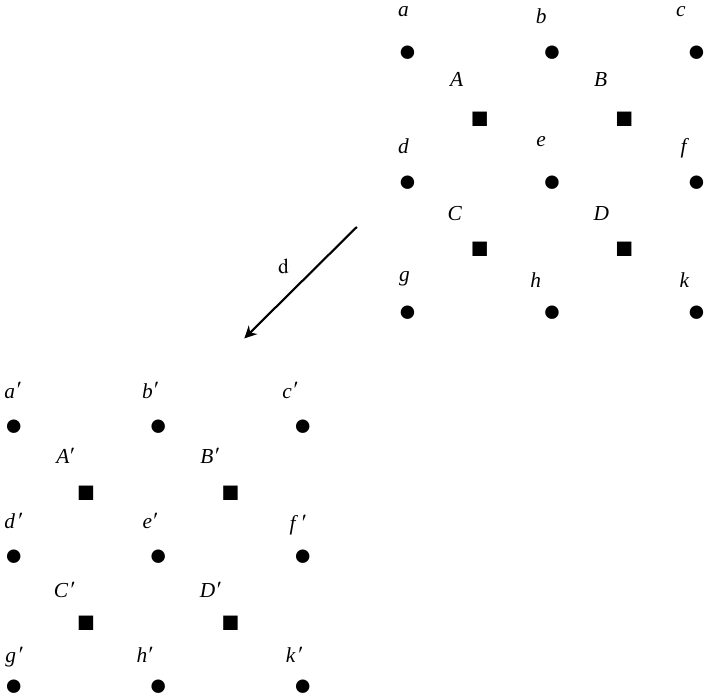}
		\caption{Pascal determinantal array $PD_k$}
		\label{fig:F6}
	\end{center}
\end{figure}

Consider the Pascal determinantal array $PD_k$, as in 
Figure~\ref{fig:F6}, where the circled-points show the entries of this array and the squared entries show the corresponding entries of the Pascal determinantal
array of order $k-1$.

Consider the previous discussions, the entries of $PD_{k+1}$ can be
obtained, as follows (see Figure~\ref{fig:F7}):

\begin{figure}[h!]
	\begin{center}
		\includegraphics[scale=0.75]{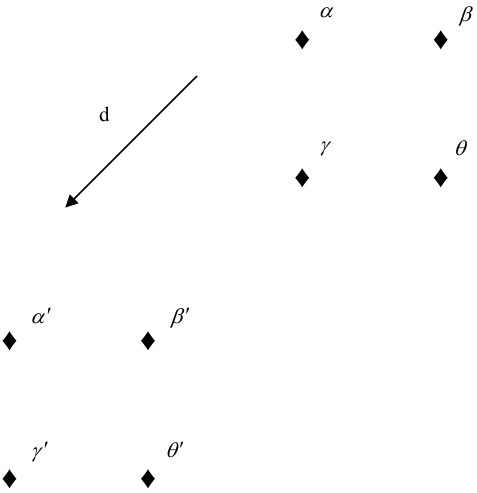}
		\caption{Pascal determinantal array $PD_{k+1}$}
		\label{fig:F7}
	\end{center}
\end{figure}

\begin{eqnarray}
\alpha &=&\frac{ae-bd}{A},~\beta=\frac{bf-ec}{B},~
\alpha' =\frac{a'e'-b'd'}{A'},~~~~\beta'=\frac{b'f'-e'c'}{B'}
\nonumber\\
\gamma &=&\frac{dh-ge}{C},~~ \theta=\frac{ek-fh}{D},~
\gamma'=\frac{d'h'-g'e'}{C'},~~ 
\theta'=\frac{e'k'-f'h'}{D'}.
\end{eqnarray}

Therefore, it is sufficient to show that
$
\frac{\alpha \theta}{\beta\gamma}
=
\frac{\alpha' \theta'}{\beta' \gamma'}
$
. 
This is equivalent to 

\begin{equation}
\frac{\frac{(ae-bd)(ek-fh)}
	{AD}}{\frac{(bf-ec)(dh-ge)}{BC}}
=
\frac{\frac{(a'e'-b'd')(e'k'-f'h')}
	{A'D'}}{\frac{(b'f'-e'c')(d'h'-g'e')}{B'C'}}
\end{equation}

By induction hypothese, we have

\[
\frac{ae}{bd}=\frac{a'e'}{b'd'},~ 
\frac{bf}{ec}=\frac{b'f'}{e'c'},~
\frac{dh}{ge}=\frac{d'h'}{g'e'},~ 
\frac{ek}{fh}=\frac{e'k'}{f'h'},~ 
\frac{AD}{BC}=\frac{A'D'}{B'C'},
\]

or equivalently

\[
\frac{ae-bd}{bd}=\frac{a'e'-b'd'}{b'd'},~
\frac{ek-fh}{fh}=\frac{e'k'-f'h'}{f'h'},
\]

\[
\frac{bf-ec}{ec}=\frac{b'f'-e'c'}{e'c'},~
\frac{dh-ge}{ge}=\frac{d'h'-g'e'}{g'e'},
\]

\[
\frac{AD}{BC}=\frac{A'D'}{B'C'}
\]

which can be simply written as

\[
\frac{\frac{(ae-bd)(ek-fh)}{(bdfh)AD}}{\frac{(bf-ec)(dh-ge)}{(ecge)BC}}
=
\frac{\frac{(a'e'-b'd')(e'k'-f'h')}{(b'd'f'h')A'D'}}{\frac{(b'f'-e'c')(d'h'-g'e')}{(e'c'g'e')B'C'}}.
\]

Therefore by $(2)$, we only need to show that

$
\frac{ecge}{bdfh}
=
\frac{e'c'g'e'}{b'd'f'h'}
$, 
or equivalently

\[
\frac{baecge}{babdfh}
=
\frac{b'a'e'c'g'e'}{b'a'b'd'f'h'}
\]

But we already know, by induction hypothesis, that

\[
\frac{ae}{bd}=\frac{a'e'}{b'd'},~~
\frac{ce}{bf}=\frac{c'e'}{b'f'},~~
\frac{bg}{ah}=\frac{b'g'}{a'h'}.
\]

Thus considering the above remark, we complete the proof by
induction.



\section{The Correctness of the Algorithm}\label{sec:correctness}

We now prove that the recursive algorithm presented in Section~\ref{sec:algorithm} correctly generates the Pascal determinantal array \(PD_k\) for every \(k \ge 1\). The proof proceeds by induction on \(k\).

\noindent \textbf{Base case (\(k=1\)):} Trivial, since \(PD_1\) is defined to be the Pascal array itself.

\noindent \textbf{Inductive step:} Assume the algorithm correctly produces \(PD_k\). That is, for all \(i,j \ge 0\),

\begin{equation}\label{Eq5}
P^{(k)}_{i,j} = \frac{P^{(k-1)}_{i+1,j+1} \; P^{(1)}_{i,j}}{P^{(k-1)}_{1,i+j+1}}.
\end{equation}

\begin{figure}[h!]
	\centering
	\includegraphics[scale=0.75]{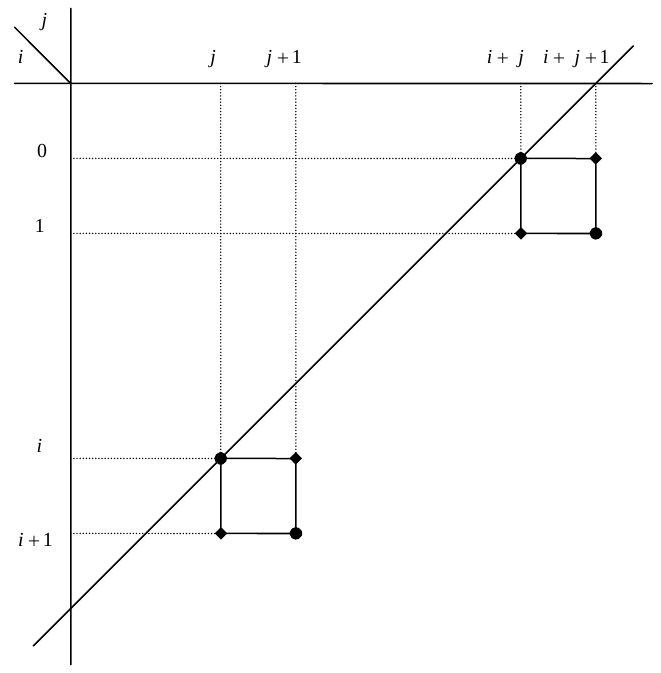}
	\caption{The sliding‑cross property applied to \(PD_k\).}
	\label{fig:F8}
\end{figure}

We must show that the same recurrence holds for \(PD_{k+1}\). Applying the sliding‑cross property (Section~\ref{sec:sliding}) to the array \(PD_k\) in the configuration shown in Figure~\ref{fig:F8} yields

\begin{equation}\label{Eq6}
\frac{P^{(k)}_{i+1,j} \; P^{(k)}_{i,j+1}}{P^{(k)}_{i,j} \; P^{(k)}_{i+1,j+1}} 
= \frac{P^{(k)}_{1,i+j} \; P^{(k)}_{0,i+j+1}}{P^{(k)}_{0,i+j} \; P^{(k)}_{1,i+j+1}}.
\end{equation}

Since \(P^{(k)}_{0,i+j+1} = P^{(k)}_{0,i+j} = 1\), equation (\ref{Eq6}) simplifies to

\begin{equation}\label{Eq7}
\frac{P^{(k)}_{1,i+j}}{P^{(k)}_{i+1,j} \; P^{(k)}_{i,j+1}}
= \frac{P^{(k-1)}_{1,i+j+1}}{P^{(k+1)}_{i,j} \; P^{(k-1)}_{i+1,j+1}}.
\end{equation}

Taking reciprocals of both sides and multiplying by \(P^{(1)}_{i,j}\) gives

\begin{equation}\label{Eq8}
\frac{P^{(1)}_{i,j} \; P^{(k)}_{i+1,j} \; P^{(k)}_{i,j+1}}{P^{(k)}_{1,i+j}}
= \frac{P^{(1)}_{i,j} \; P^{(k+1)}_{i,j} \; P^{(k-1)}_{i+1,j+1}}{P^{(k-1)}_{1,i+j+1}}.
\end{equation}

Using the induction hypothesis (\ref{Eq5}) on the right‑hand side, we obtain

\begin{equation}\label{Eq9}
\frac{P^{(1)}_{i,j} \; P^{(k)}_{i+1,j} \; P^{(k)}_{i,j+1}}{P^{(k)}_{1,i+j}}
= P^{(k+1)}_{i,j} \; P^{(k)}_{i,j}.
\end{equation}

Divide both sides of (\ref{Eq9}) by \(P^{(k)}_{i,j}\) and multiply the left-hand side by \(\frac{P^{(k)}_{i+1,j+1}}{P^{(k)}_{i+1,j+1}}\) (which equals 1):

\begin{equation}\label{Eq10}
\frac{P^{(k)}_{i+1,j+1} \; P^{(1)}_{i,j} \; P^{(k)}_{i+1,j} \; P^{(k)}_{i,j+1}}
{P^{(k)}_{i+1,j+1} \; P^{(k)}_{i,j} \; P^{(k)}_{1,i+j}}
= P^{(k+1)}_{i,j}.
\end{equation}

Now apply (\ref{Eq5}) again to replace the product \(P^{(k)}_{i+1,j} \; P^{(k)}_{i,j+1}\) in the numerator:

\begin{equation}\label{Eq11}
\frac{P^{(k)}_{i+1,j+1} \; P^{(1)}_{i,j} \; P^{(k)}_{1,i+j} \; P^{(k)}_{0,i+j+1}}
{P^{(k+1)}_{0,i+j} \; P^{(k)}_{1,i+j+1} \; P^{(k)}_{1,i+j}}
= P^{(k+1)}_{i,j}.
\end{equation}

Recall that \(P^{(k)}_{0,i+j+1} = P^{(k+1)}_{0,i+j} = 1\). Cancelling the common factor \(P^{(k)}_{1,i+j}\) from numerator and denominator yields

\[
\frac{P^{(k)}_{i+1,j+1} \; P^{(1)}_{i,j}}{P^{(k)}_{1,i+j+1}} = P^{(k+1)}_{i,j}.
\]

This is precisely the recurrence required by the algorithm for producing \(PD_{k+1}\) from \(PD_k\). Hence the algorithm is correct for all \(k \ge 1\).



\section{A Geometric Interpretation of $P^{(k)}_{i,j}$}\label{sec:geometry}

In this section we present a geometric interpretation of the $(i,j)$-entry of the Pascal determinantal array of order $k$. We first derive a new identity for Pascal array entries as a consequence of the sliding-cross rule (Section~\ref{sec:sliding}).

\begin{figure}[h!]
	\centering
	\includegraphics[scale=0.8]{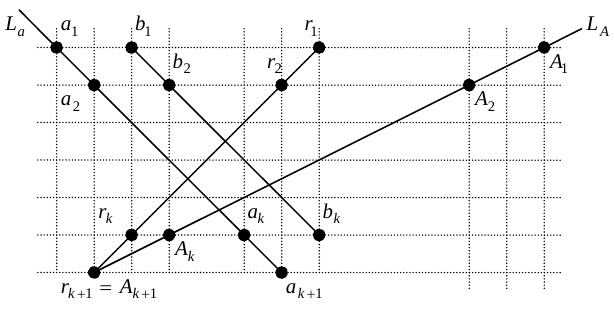}
	\caption{Geometric interpretation of the new identity.}
	\label{fig:F9}
\end{figure}

For every positive integer $k$ and nonnegative integer $j$, we have
\[
P_{0,j+2k+1} \cdot P_{1,j+2k-1} \cdot \ldots \cdot P_{k-1,j+3} \cdot P_{k,j+1}
= P_{0,j} \cdot P_{1,j+1} \cdot \ldots \cdot P_{k-1,j+k-1} \cdot P_{k,j+k}.
\]
The corresponding entries are shown in Figure~\ref{fig:F9}. For simplicity, denote the entries on the left-hand side (LHS) by $a_i$ ($1 \le i \le k+1$) and those on the right-hand side (RHS) by $A_i$ ($1 \le i \le k+1$). Note that $a_1 = A_1 = 1$. The LHS entries lie on a line $L_a$ with slope $-\frac{1}{1}$, while the RHS entries lie on a line $L_A$ with slope $+\frac{1}{2}$.

We give a geometric proof of this identity using induction on $k$ and the sliding-cross rule. The base case ($k=1$) is trivial. Assume the identity holds for entries $b_1,b_2,\ldots,b_k$ and $A_1,A_2,\ldots,A_k$, i.e.,
\begin{equation}\label{keycross1}
b_1 b_2 \ldots b_k = A_1 A_2 \ldots A_k.
\end{equation}

Consider two parallel crosses $(b_1,\ldots,b_k,r_1,\ldots,r_k)$ and $(a_2,\ldots,a_{k+1},r_2,\ldots,r_{k+1})$ of size $k$. By the sliding-cross rule,
\[
\frac{b_1 b_2 \ldots b_k}{r_1 r_2 \ldots r_k} = \frac{a_2 a_3 \ldots a_{k+1}}{r_2 r_3 \ldots r_{k+1}},
\]
or equivalently,
\begin{equation}\label{keycross2}
b_1 b_2 \ldots b_k \cdot A_{k+1} = a_1 a_2 \ldots a_{k+1}.
\end{equation}

Combining (\ref{keycross1}) and (\ref{keycross2}) yields
\[
a_1 a_2 \ldots a_{k+1} = A_1 A_2 \ldots A_{k+1},
\]
completing the induction.

\subsection*{Closed‑form Formula for $P^{(k)}_{i,j}$}

Using the recursive algorithm from Section~\ref{sec:algorithm},
\[
P^{(k+1)}_{i,j} = \frac{P^{(k)}_{i+1,j+1} \; P^{(1)}_{i,j}}{P^{(k)}_{1,i+j+1}}.
\]

Iterating this recurrence $k-1$ times gives
\[
P^{(k+1)}_{i,j} = 
\frac{
	P^{(1)}_{i+k,j+k} \cdot P^{(1)}_{i+(k-1),j+(k-1)} \cdot \ldots \cdot P^{(1)}_{i+1,j+1} \cdot P^{(1)}_{i,j}
}{
	P^{(1)}_{1,i+j+2k-1} \cdot 
	P^{(2)}_{1,i+j+2k-3} \cdot \ldots \cdot 
	P^{(k-1)}_{1,i+j+3} \cdot 
	P^{(k)}_{1,i+j+1}
}.
\]

By Rahimpour's theorem \cite{m1}, $P^{(k)}_{1,j} = P^{(1)}_{k,j}$ for all $j \ge 0$ and $k \ge 1$. Hence
\[
P^{(k+1)}_{i,j} = 
\frac{
	P_{i+k,j+k} \cdot P_{i+(k-1),j+(k-1)} \cdot \ldots \cdot P_{i+1,j+1} \cdot P_{i,j}
}{
	P_{1,i+j+2k-1} \cdot P_{2,i+j+2k-3} \cdot \ldots \cdot P_{k-1,i+j+3} \cdot P_{k,i+j+1}
}.
\]

Applying the new identity with $j$ replaced by $i+j$ yields an alternative form:
\[
P^{(k+1)}_{i,j} = 
\frac{
	P_{i+k,j+k} \cdot P_{i+(k-1),j+(k-1)} \cdot \ldots \cdot P_{i+1,j+1} \cdot P_{i,j}
}{
	P_{1,i+j+1} \cdot P_{2,i+j+2} \cdot \ldots \cdot P_{k-1,i+j+(k-1)} \cdot P_{k,i+j+k}
}.
\]

Now consider two parallel symmetric crosses of size $k+1$:
\[
(P_{i,j},\ldots,P_{i+k,j+k},P_{i,j+k},\ldots,P_{i+k,j})
\]
and
\[
(P_{0,i+j},\ldots,P_{k,i+j+k},P_{k,i+j},\ldots,P_{0,i+j+k}).
\]
Using the sliding property, we obtain
\[
P^{(k+1)}_{i,j} = 
\frac{
	P_{i+k,j} \cdot P_{i+(k-1),j+1} \cdot \ldots \cdot P_{i+1,j+(k-1)} \cdot P_{i,j+k}
}{
	P_{0,i+j+k} \cdot P_{1,i+j+(k-1)} \cdot \ldots \cdot P_{k-1,i+j+1} \cdot P_{k,i+j}
}.
\]

Finally, replacing $k$ by $k-1$ gives the closed‑form expression for $P^{(k)}_{i,j}$:
\begin{equation}\label{eq:final-form}
P^{(k)}_{i,j} = 
\frac{
	P_{i+(k-1),j} \cdot P_{i+(k-2),j+1} \cdot \ldots \cdot P_{i+1,j+(k-2)} \cdot P_{i,j+(k-1)}
}{
	P_{i+j+(k-1),0} \cdot P_{i+j+(k-2),1} \cdot \ldots \cdot P_{i+j+1,k-2} \cdot P_{i+j,k-1}
}.
\end{equation}

\subsection*{Double‑stick Interpretation}

\begin{figure}[h!]
	\centering
	\includegraphics[scale=0.75]{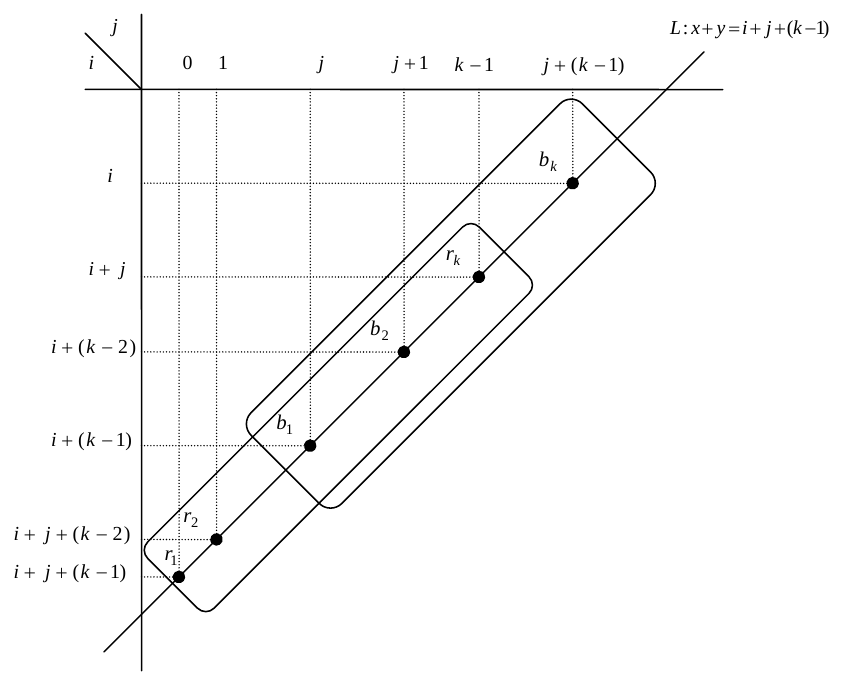}
	\caption{Double‑stick interpretation of $P^{(k)}_{i,j}$.}
	\label{fig:F10}
\end{figure}

Equation (\ref{eq:final-form}) admits a vivid geometric interpretation. Define a \emph{double stick} (Figure~\ref{fig:F10}) as a pair of sequences
\[
(b_1,b_2,\ldots,b_k \mid r_1,r_2,\ldots,r_k),
\]
where all entries lie on the discrete line $L: x+y = i+j+(k-1)$, with
\[
b_1 = P_{i+(k-1),j},\; b_k = P_{i,j+(k-1)},\qquad
r_1 = P_{i+j+(k-1),0},\; r_k = P_{i+j,k-1}.
\]

The \emph{weight} of this double stick is
\[
W^{(k)}_{i,j} = \frac{b_1 b_2 \ldots b_k}{r_1 r_2 \ldots r_k}.
\]

Comparing with (\ref{eq:final-form}) shows that $P^{(k)}_{i,j} = W^{(k)}_{i,j}$. Thus the determinant $P^{(k)}_{i,j}$ equals the weight of the corresponding double stick. The fact that the overlapping area of the double stick grows with $k$ underlies the proof of our main conjecture.

\subsection*{Geometric Proof of the Conjecture}

In the language of Pascal determinantal arrays, Conjecture~\ref{HasanConj} states that for any $i,k,n$ with $n \ge k$,
\begin{equation}\label{eq:conj-statement}
P^{(n)}_{i,k} = P^{(k)}_{i,n}.
\end{equation}

\begin{figure}[h!]
	\centering
	\includegraphics[scale=.8]{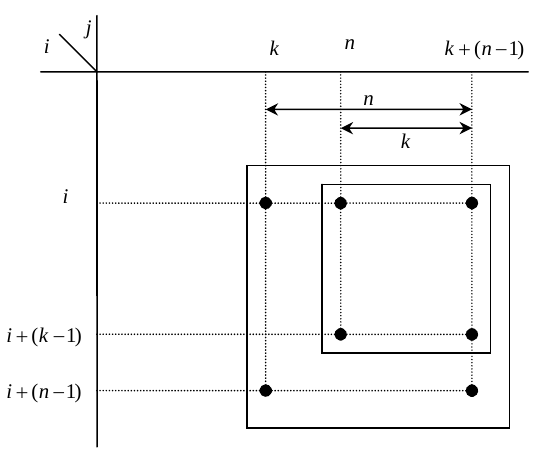}
	\caption{Sub‑arrays appearing in the conjecture.}
	\label{fig:F11}
\end{figure}

\begin{figure}[h!]
	\centering
	\includegraphics[scale=0.7]{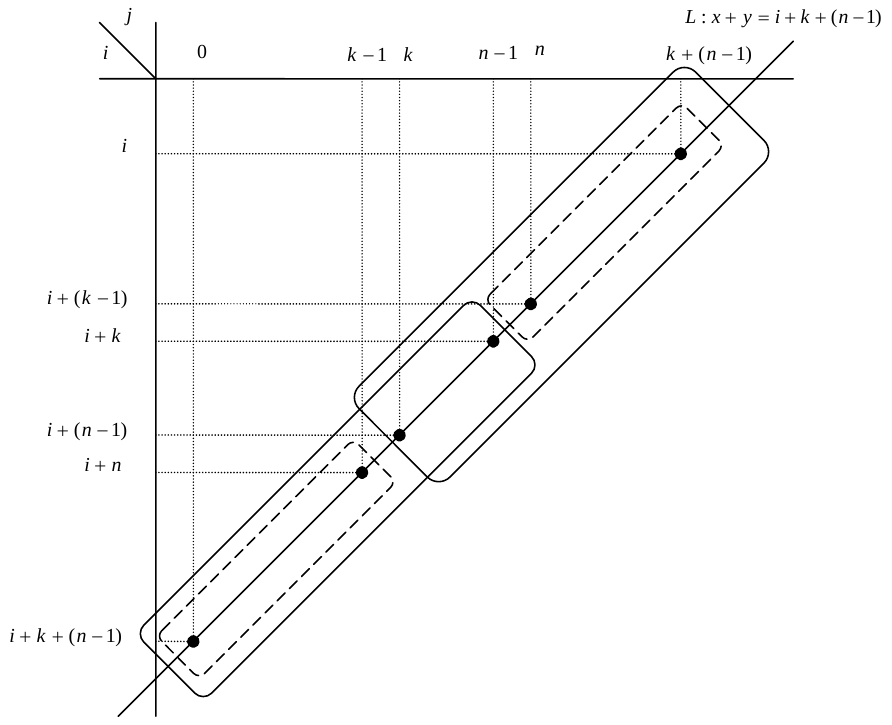}
	\caption{Geometric proof of the conjecture: cancelling the overlapping region (solid) leaves the dashed double stick.}
	\label{fig:F12}
\end{figure}

Geometrically, $P^{(n)}_{i,k}$ is the weight $W^{(n)}_{i,k}$ of the solid double stick in Figure~\ref{fig:F12}, which lies on the line $L: x+y = i+k+(n-1)$. Explicitly,
\[
P^{(n)}_{i,k} = W^{(n)}_{i,k} = 
\frac{
	P_{i+(n-1),k} \cdot \ldots \cdot P_{i,k+(n-1)}
}{
	P_{i+k+(n-1),0} \cdot \ldots \cdot P_{i+k,n-1}
}.
\]

Equivalently, we may write
\[
P^{(n)}_{i,k} = 
\frac{
	\big(P_{i+(n-1),k} \ldots P_{i+k,n-1}\big) \;
	\big(P_{i+(k-1),n} \ldots P_{i,n+(k-1)}\big)
}{
	\big(P_{i+n+(k-1),0} \ldots P_{i+n,k-1}\big) \;
	\big(P_{i+(n-1),k} \ldots P_{i+k,n-1}\big)
}.
\]

Cancelling the common factor $\big(P_{i+(n-1),k} \ldots P_{i+k,n-1}\big)$ (which corresponds geometrically to removing the overlapping area of the original double stick) yields
\[
P^{(n)}_{i,k} = 
\frac{
	P_{i+(k-1),n} \cdot \ldots \cdot P_{i,n+(k-1)}
}{
	P_{i+n+(k-1),0} \cdot \ldots \cdot P_{i+n,k-1}
}.
\]

The right‑hand side is precisely the weight $W^{(k)}_{i,n}$ of the dashed double stick in Figure~\ref{fig:F12}, i.e., $P^{(k)}_{i,n}$. Hence
\[
P^{(n)}_{i,k} = P^{(k)}_{i,n},
\]
which completes the geometric proof of Conjecture~\ref{HasanConj}.

\newpage


\section{Conclusion and Connections to Prior Work}

In this paper, we introduced Pascal determinantal arrays $PD_k$, a new infinite family of arrays derived from $k \times k$ minors of the classical Pascal array. We provided a recursive algorithm for generating these arrays, proved its correctness using Dodgson's condensation and a weighted sliding-cross rule, and gave a geometric interpretation of $P^{(k)}_{i,j}$ as the weight of a \emph{double stick} in the Pascal plane. Our main result established the symmetric identity $P^{(k)}_{i,j} = P^{(j)}_{i,k}$, which generalizes Rahimpour's original determinantal conjecture.

\subsection*{Algebraic Connections and Alternative Proofs}

After the initial submission of this work to \emph{ArXiv}, we learned from colleagues of important connections to previously studied combinatorial objects. 
Christian Krattenthaler \cite{kratt} observed that $P^{(k)}_{i,j}$ admits a closed-form evaluation:
\begin{equation}\label{eq:kratt-formula}
P^{(k)}_{i,j} = \prod_{\ell=0}^{k-1} \frac{(i+j+\ell)! \, \ell!}{(i+\ell)! \, (j+\ell)!},
\end{equation}
which follows from standard determinant techniques (Pochhammer symbols and Vandermonde determinants). This formula yields an immediate algebraic proof of our main identity $P^{(k)}_{i,j} = P^{(j)}_{i,k}$ via induction.

Johann Cigler \cite{cigler} further noted that \eqref{eq:kratt-formula} coincides exactly with the \emph{Hoggatt binomials} $\binom{i+j}{j}_k$ introduced by Fielder and Alford in 1989. In their notation,
\[
P^{(k)}_{i,j} = \binom{i+j}{j}_k = \frac{\langle i+j\rangle_k!}{\langle i\rangle_k! \, \langle j\rangle_k!},
\]
where $\langle n\rangle_k = \prod_{i=1}^n \binom{i+k-1}{k}$. Many of the properties we derived geometrically---such as the sliding cross rule---are naturally encoded in the algebraic identities satisfied by Hoggatt binomials.

\subsection*{Contributions of This Work}

While the algebraic expression \eqref{eq:kratt-formula} was known in other contexts, our contributions are distinct and complementary:

\begin{itemize}
	\item We introduced the \emph{determinantal definition} $P^{(k)}_{i,j} = \det\big(\binom{i+j+\alpha+\beta}{i+\alpha}\big)_{0\le\alpha,\beta\le k-1}$, which provides a natural generalization of the Pascal array and reveals its layered minor structure.
	
	\item We developed a \emph{geometric framework} based on weighted double sticks and sliding crosses, offering visual intuition and pictorial proofs that are independent of algebraic computation.
	
	\item We established \emph{recursive generation algorithms} and proved their correctness using Dodgson condensation and the sliding property, providing constructive tools for working with these arrays.
	
	\item We demonstrated how classical combinatorial results---Dodgson's condensation, the star of David rule, and Vandermonde determinants---interact within this geometric setting.
\end{itemize}

Thus, our work provides a \emph{geometric and algorithmic} perspective on structures that were previously studied primarily through algebraic and enumerative methods.

\subsection*{Future Directions}

Several natural avenues for further investigation arise from our results:

\begin{enumerate}
	\item \textbf{Generalized Pascal arrays:} Replace the binomial coefficients $\binom{i+j}{i}$ with other two-parameter sequences (e.g., $q$-binomials, Fibonomial coefficients) and study the corresponding determinantal arrays. Do they admit similar geometric interpretations or sliding properties?
	
	\item \textbf{Higher-dimensional analogs:} Extend the double-stick construction to higher dimensions, perhaps using hypercubes or simplices in multidimensional Pascal pyramids.
	
	\item \textbf{Combinatorial interpretations:} Find explicit combinatorial models (lattice paths, tableaux, etc.) for the entries of $PD_k$ when $k>2$, similar to the Narayana interpretation for $PD_2$.
	
	\item \textbf{Algebraic geometry connections:} Explore whether the determinantal identities have interpretations in terms of Schubert calculus, Grassmannians, or determinantal varieties.
	
	\item \textbf{Algorithmic applications:} Investigate computational aspects---efficient computation of large determinants in $PD_k$, connections to iterative matrix methods, or potential applications in coding theory.
\end{enumerate}

The interplay between algebraic closed forms, geometric visualization, and recursive algorithms exhibited in this paper suggests that Pascal determinantal arrays offer a rich testbed for exploring deeper connections between combinatorics, algebra, and geometry.






\begin{thebibliography}{99}
	\bibitem{m1}
	H. Teimoori and M. Bayat, Fermat Row-Eliminated Matrices and Some Binomial Determinants, {\it J. Math.
		Gazette.}, 87 (2003), pp. 114-118.
	\bibitem{m2}
	M. Bayat and H. Teimoori, Pascal k-eliminated Functional Matrix and It's Property, {\it Linear Algebra Appl}, 308 (2000), pp.65-75.
	\bibitem{m3}
	P. Hilton and J. Pedersen, Relation Geometry and Algebra in the Pascal Triangle, Hexagon, Tetrahedron, and Cuboctahedron Part I: Binomial Coefficients, Extended Binomial Cofficients and Prepration for further Work,{\it The College Math. J}, 30 (3)(1999), pp. 170-186.
	\bibitem{m4}
	C. L. Dodgson, Condensation of Determinants: Being a New and Brief
	Method for Computing their Arithmetic Values, {\it Proc. Roy. Soc.
		Ser. A}, 15 (1866)150-155.
	\bibitem{m5}
	N.J.A. Sloane, The On-Line Encyclopedia of Integer Sequences,
	\url{https://oeis.org}, 2008.
	\bibitem{m6}
	K. Robin McLean, Determinants of binomial coefficients, {\it J.
		Math. Gazette.}, 89 (2005), pp. 247-250.
	
	\bibitem{kratt}
	C. Krattenthaler, \emph{Advanced determinant calculus}, S\'eminaire Lotharingien Combin. 42 (``The Andrews Festschrift'') (1999), Article B42q, 67~pp.
	
	\bibitem{fielder-alford}
	D. C. Fielder and C. O. Alford, \emph{On a conjecture by Hoggatt with extensions to Hoggatt sums and Hoggatt triangles}, Fibonacci Quart. 27 (1989), no.~2, 160--168.
	
	\bibitem{cigler}
	J. Cigler, Pascal triangle, Hoggatt matrices, and analogous constructions, (2021), \url{arXiv:2103.01652v3}
	
	
\end{thebibliography}
\end{document}